\documentclass[10pt]{amsart}

\usepackage{amsmath,xspace,amssymb,mathrsfs}
\usepackage{color}

\input xy
\xyoption{all}
\xyoption{2cell}
\UseAllTwocells
\CompileMatrices

\newcommand{\Spec}{\operatorname{Spec}}
\renewcommand{\phi}{\varphi}

\newcommand{\Ker}{\operatorname{Ker}}

\newcommand{\Supp}{\operatorname{Supp}}

\newcommand{\Fun}{\operatorname{Fun}}
\newcommand{\D}{\operatorname{D}}

\newcommand{\Fin}{\operatorname{Fin}}

\newcommand{\Su}{\operatorname{S}}

\newtheorem{proposition}{Proposition}[section]
\newtheorem{lemma}[proposition]{Lemma}

\newtheorem{corollary}[proposition]{Corollary}
\newtheorem{theorem}[proposition]{Theorem}

\theoremstyle{definition}

\newtheorem{remark}[proposition]{Remark}

\usepackage{etoolbox}
\makeatletter
\patchcmd{\@settitle}{\uppercasenonmath\@title}{}{}{}
\patchcmd{\@setauthors}{\MakeUppercase}{}{}{}
\makeatother

\begin{document}

\title[On the direct product of fields]{On the direct product of fields with an application}

\author[A. Tarizadeh]{Abolfazl Tarizadeh}
\address{Department of Mathematics, Faculty of Basic Sciences, University of Maragheh \\
P. O. Box 55136-553, Maragheh, Iran.
 }
\email{ebulfez1978@gmail.com \\ atarizadeh@maragheh.ac.ir}

\date{}
\subjclass[2010]{14A05, 14A15, 13A15, 13B30}
\keywords{Direct product of fields; Separated scheme}

\begin{abstract} In this paper, the (infinite) direct product of fields is investigated. In particular, the finiteness of a given set is characterized in terms of some ring-theoretic observations. Next, a certain localization (whose multiplicative set formed by cofinite sets) of the direct product of fields is studied. Finally, it is shown that every set $X$ can be made into a separated scheme, and this scheme is an affine scheme if and only if $X$ is a finite set.
\end{abstract}

\maketitle

\section{Introduction}

The structure of prime ideals of an infinite direct product of commutative rings (even fields) is amazingly complicated and includes a huge number of prime ideals. In the literature, these rings have been investigated and several interesting results have been obtained. For further studies on these topics see the References (\cite{Anderson}, \cite{Gilmer-Heinzer}, \cite{Gilmer-Heinzer 2}, \cite{Gilmer-Heinzer 3}, \cite{Maroscia}, \cite{Ronnie et al} and \cite{Tarizadeh 3}). But there is a pressing need for new constructions and approaches in order to understand the structure of prime ideals of infinite direct product rings more deeply. There are some hopes in order to realize these goals. For instance, the space $\Spec(\Lambda)$ is the Stone-$\check{C}ech$ compactification of the discrete space $X$ where the ring $\Lambda$ is the direct product of a family fields indexed by the set $X$. (It is worth mentioning that in \cite[Theorems 3.5 and 5.4]{Tarizadeh} we improve this fundamental result by passing from fields to integral domains and local rings). This important result together with the results of  \cite{Abolfazl-Taheri} are the main motivations in writing the present paper. The other motivation comes from a geometric point of view. In fact in this paper, after studying some basic properties of the ring $\Lambda$, then as an application, we assign to every set $X$ a separated scheme structure whose structure sheaf is obtained by the direct product of fields. It is also shown that this scheme is an affine scheme if and only if $X$ is a finite set (see Theorem \ref{Theorem set-scheme}). In Theorem \ref{Proposition DI}, finiteness of a given set is characterized in terms of some ring-theoretic notions. In Theorem \ref{Theorem IVD2}, it is proved that the localization $T^{-1}\Lambda$ is canonically isomorphic to a certain quotient of $\Lambda$ where the multiplicative set $T$ formed by cofinite sets. Theorems \ref{Proposition DIII} and \ref{Th Del} are further main results of this paper. In a subsequent work \cite{Tarizadeh 3}, we study various aspects, especially the structure of prime ideals of infinite direct products of general commutative rings.

\section{Preliminaries}

We collect in this section some basic background for the reader's convenience. \\

Throughout the paper, $\Lambda=\prod\limits_{x\in X}K_{x}$ where $X$ is an index set and each $K_{x}$ is a field. For each $f=(f_{x})\in\Lambda$, the subset $\Supp(f)=\{x\in X: f_{x}\neq0\}$ is simply denoted by $\Su(f)$. \\

If $X$ is a set then its power set $\mathcal{P}(X)$ together with the symmetric difference $A+B=(A\cup B)\setminus (A\cap B)$ as the addition and the intersection $A.B=A\cap B$ as the multiplication forms a commutative ring whose zero and unit are respectively the empty set and the whole set $X$. The ring $\mathcal{P}(X)$ is called the \emph{power set ring} of $X$. If $f:X\rightarrow Y$ is a function then the map $\mathcal{P}(f):\mathcal{P}(Y)\rightarrow\mathcal{P}(X)$ defined by $A\mapsto f^{-1}(A)$ is a morphism of rings. In fact, the assignments $X\mapsto\mathcal{P}(X)$ and $f\mapsto\mathcal{P}(f)$ form a faithful contravariant functor from the category of sets to the category of Boolean rings. We call it the \emph{power set functor}. By $\Fin(X)$ we mean the set of all finite subsets of $X$. It is an ideal of $\mathcal{P}(X)$. \\

The set of all functions from a set $X$ to a nonzero ring $K$, denoted $\Fun(X,K)$, with the usual addition and multiplication of functions is a commutative ring. This ring is canonically isomorphic to $\prod\limits_{x\in X}K$.
If $f:X\rightarrow Y$ is a function then the induced map $\Fun(f):\Fun(Y,K)\rightarrow\Fun(X,K)$ given by $g\mapsto g\circ f$ is a morphism of rings. Indeed, $\Fun(-,K)$ is a faithful contravariant functor from the category of sets to category of commutative rings. To see its faithfulness, suppose $\Fun(f)=\Fun(g)$ for some functions $f,g:X\rightarrow Y$. We show that $f=g$. For each $x\in X$ consider the function $h_{x}:Y\rightarrow K$ with $h_{x}(y)=\delta_{f(x),y}$ where $\delta_{f(x),y}$ is the Kronecker delta. We have $h_{x}\circ f=h_{x}\circ g$. In particular, $1=\delta_{f(x),f(x)}=\delta_{f(x),g(x)}$. This yields that $f(x)=g(x)$. \\

The map $\mathcal{P}(X)\rightarrow\Fun(X,\mathbb{Z}_{2})$ given by $A\mapsto\chi_{A}$ is an isomorphism of rings where $\chi_{A}$ is the characteristic function of $A$ and $\mathbb{Z}_{2}=\{0,1\}$ is the field of integers modulo 2. \\

If $f$ is a member of a commutative ring $R$ then $D(f)=\{\mathfrak{p}\in\Spec(R): f\notin\mathfrak{p}\}$ and by $(f)$ we mean the ideal of $R$ generated by $f$. \\

It is well known that a scheme $(X,\mathscr{O}_{X})$ is a separated scheme if and only if for every affine opens $U$ and $V$ of $X$ then $U\cap V$ is an affine open and the canonical ring map $\mathscr{O}_{X}(U)\otimes_{\mathscr{O}_{X}(X)}\mathscr{O}_{X}(V)
\rightarrow\mathscr{O}_{X}(U\cap V)$ which sends each pure tensor $s\otimes t$ into $(s|_{U\cap V})\cdot(t|_{U\cap V})$ is surjective where $\mathscr{O}_{X}(X)$ is the ring of global sections.

\section{Some observations on the ring $\Lambda$}

First note that the map $\Lambda=\prod\limits_{x\in X}K_{x}\rightarrow\mathcal{P}(X)$ given by $f\mapsto\Su(f)$ is multiplicative: $\Su(fg)=\Su(f)\cap\Su(g)$ for all $f,g\in\Lambda$. But in general it is not additive. In fact, $\Su(f)+\Su(g)\subseteq\Su(f+g)\subseteq\Su(f)\cup\Su(g)$.
We have also the following result.

\begin{theorem}\label{Proposition DIII} If $f,g\in\Lambda$ then the following assertions hold. \\
$\mathbf{(i)}$ $(f)=\{h\in\Lambda: \Su(h)\subseteq\Su(f)\}$ and the ring $\Lambda/(f)$ is canonically isomorphic to $\prod\limits_{x\in X\setminus\Su(f)}K_{x}$. \\
$\mathbf{(ii)}$ $(f)=(g)$ if and only if $\Su(f)=\Su(g)$. \\
$\mathbf{(iii)}$ Every ideal of $\Lambda$ is a radical ideal. \\
$\mathbf{(iv)}$ $\D(f)\subseteq\D(g)$ if and only if $\Su(f)\subseteq\Su(g)$.
\end{theorem}

\begin{proof} (i): If $\Su(h)\subseteq\Su(f)$ then $h=h'f$ where $h'_{x}$ is either $h_{x}f^{-1}_{x}$ or $0$, according as $x\in\Su(h)$ or $x\notin\Su(h)$. To see the isomorphism consider the canonical projection $\phi:\Lambda\rightarrow\prod\limits_{x\in X\setminus\Su(f)}K_{x}$ given by $h\mapsto(h_{x})_{x\in X\setminus\Su(f)}$. This ring map is surjective and $\Ker\phi=(f)$. \\ 
(ii): It follows from (i). \\
(iii): Let $I$ be an ideal of $\Lambda$. If $f\in\sqrt{I}$ then there exists some $n\geq1$ such that $f^{n}\in I$. But $\Su(f^{n})=\Su(f)$. Thus by (ii), $f\in(f^{n})\subseteq I$. \\
(iv): It follows from the fact that $\D(f)\subseteq\D(g)$ if and only if $f\in\sqrt{(g)}$.
\end{proof}

It is obvious that an element $f\in\Lambda$ is invertible (in $\Lambda$)
if and only if $\Su(f)=X$. \\

Let $\mathfrak{F}(\Lambda)$ be the set of all $f\in\Lambda$ such that $\Su(f)$ is a finite set. Then $\mathfrak{F}(\Lambda)$ is an ideal of $\Lambda$ and it is generated by the sequences $\Delta_{x}=(\delta_{x,y})_{y\in X}$ where $x\in X$ and $\delta_{x,y}$ is the Kronecker delta. Indeed, if $f\in\mathfrak{F}(\Lambda)$ then we may write $f=\sum\limits_{x\in\Su(f)}r_{x}\Delta_{x}$ where $r_{x}$ is the sequence in $\Lambda$ whose $x$-th component is $f_{x}$ and the other components are zero.

\begin{proposition}\label{Proposition DII} Let $R=\prod\limits_{x\in\Su(f)}K_{x}$ and $R'=\prod\limits_{x\in\Su(g)}K_{x}$ where $f,g\in\Lambda$. Then the following assertions hold. \\
$\mathbf{(i)}$ If $\Su(f)\cap\Su(g)=\emptyset$ then the ring $\prod\limits_{x\in\Su(f)\cup\Su(g)}K_{x}$ is canonically isomorphic to $R\times R'$. \\
$\mathbf{(ii)}$ The $\Lambda-$algebras $R\otimes_{\Lambda}R'$ and $\prod\limits_{x\in\Su(fg)}K_{x}$ are canonically isomorphic. \\
$\mathbf{(iii)}$ The $\Lambda-$algebras $\Lambda/\mathfrak{F}(\Lambda)\otimes_{\Lambda}R$ and $R/\mathfrak{F}(R)$ are canonically isomorphic.
\end{proposition}

\begin{proof} (i): It is clear. \\
(ii): Consider the sequence $h\in\Lambda$ where $h_{x}$ is either $0$ or $1$, according as $x\in\Su(g)$ or $x\notin\Su(g)$. Then by Proposition \ref{Proposition DIII}(i), the ring $R'$ is isomorphic to $\Lambda/(h)$. Then we have the following isomorphisms of rings: $$R\otimes_{\Lambda}R'\simeq R/(h)^{e}$$ where $(h)^{e}$ is the extension of the ideal $(h)$ under the canonical projection $\Lambda\rightarrow R$. But the canonical projection $R\rightarrow\prod\limits_{x\in\Su(fg)}K_{x}$ is surjective and whose kernel is equal to $(h)^{e}$. \\
(iii): It is proved exactly like (ii).
\end{proof}

The above theorem, in particular, yields the following special case which deserves to be mentioned separately.

\begin{corollary} If $A,B\in R=\mathcal{P}(X)$ then the following assertions hold. \\
$\mathbf{(i)}$ If $A\cap B=\emptyset$ then the ring $\mathcal{P}(A\cup B)$ is canonically isomorphic to $\mathcal{P}(A)\times\mathcal{P}(B)$. In particular,  $R\simeq\mathcal{P}(A)\times\mathcal{P}(A^{c})$. \\
$\mathbf{(ii)}$ The ring $\mathcal{P}(A+B)$ is canonically isomorphic to $\mathcal{P}(A\setminus B)\times\mathcal{P}(B\setminus A)$. \\
$\mathbf{(iii)}$ The $R-$algebras $\mathcal{P}(A\cap B)$ and $\mathcal{P}(A)\otimes_{R}\mathcal{P}(B)$ are canonically isomorphic. In particular, $\mathcal{P}(A)\otimes_{R}\mathcal{P}(A^{c})=0$. \\
$\mathbf{(iv)}$ The $R-$algebras $R/\Fin(X)\otimes_{R}\mathcal{P}(A)$ and $\mathcal{P}(A)/\Fin(A)$ are canonically isomorphic. In particular, if $A$ is a finite set then $R/\Fin(X)\otimes_{R}\mathcal{P}(A)=0$.
\end{corollary}

It is easy to see that every prime ideal of $\Lambda$ is a maximal ideal (i.e., it is a zero dimensional ring). It follows that for each $x\in X$ then $\mathfrak{m}_{x}=\{f\in\Lambda: f_{x}=0\}$ is a maximal ideal of $\Lambda$ and it is generated by the sequence $1-\Delta_{x}$. Indeed, $f=f(1-\Delta_{x})$ for all $f\in\mathfrak{m}_{x}$.

\begin{theorem}\label{Proposition DI} For a set $X$ the following statements are equivalent. \\
$\mathbf{(i)}$ $X$ is a finite set. \\
$\mathbf{(ii)}$ The zero ideal of $\Lambda$ is a finite intersection of maximal ideals of $\Lambda$. \\
$\mathbf{(iii)}$ The maximal ideals of $\Lambda$ are precisely of the form $\mathfrak{m}_{x}$ with $x\in X$. \\
$\mathbf{(iv)}$ $\Lambda$ is a Noetherian ring. \\
$\mathbf{(v)}$ $\Lambda$ is an Artinian ring.
\end{theorem}

\begin{proof} (i)$\Rightarrow$(ii): For any set $X$ we have
$\bigcap\limits_{x\in X}\mathfrak{m}_{x}=0$. \\
(ii)$\Rightarrow$(i): By the hypothesis there exists a finite set $\{M_{1},...,M_{n}\}$ of (distinct) maximal ideals of $\Lambda$ such that $\bigcap\limits_{k=1}^{n}M_{k}=0$. We shall prove that $X$ has (at most) $n$ elements. Choose some $x_{1}\in X$ then there exists some $k$, say $k=1$, such that $M_{1}\subseteq\mathfrak{m}_{x_{1}}$. This yields that $M_{1}=\mathfrak{m}_{x_{1}}$. If $X\setminus\{x_{1}\}$ is the empty set then we are done, otherwise we may choose some $x_{2}\in X\setminus\{x_{1}\}$. Then there exists some $k\in\{2,...,n\}$, say $k=2$, such that $M_{2}=\mathfrak{m}_{x_{2}}$. Hence, in this way we may find (distinct) elements $x_{1},...,x_{n}\in X$ such that $M_{k}=\mathfrak{m}_{x_{k}}$ for all $k$. If $X\setminus\{x_{1},...,x_{n}\}$ is non-empty then choose some $x$ in it. Similarly as the above, there exists some $k$ such that $\mathfrak{m}_{x}=\mathfrak{m}_{x_{k}}$. Therefore $x=x_{k}$ which is a contradiction. Thus $X$ has (at most) $n$ elements. \\
(i)$\Rightarrow$(iii): Let $M$ be a maximal ideal of $\Lambda$. Since $X$ is a finite set, we may write $1_{R}=\sum\limits_{x\in X}\Delta_{x}$.  But $1_{R}\notin M$. Hence there exists some $x\in X$ such that $\Delta_{x}\notin M$. It follows that $1-\Delta_{x}\in M$ and so $\mathfrak{m}_{x}=M$. \\
(iii)$\Rightarrow$(iv): By the hypothesis, every maximal ideal of $R$ is of the form $\mathfrak{m}_{x}$ which is a principal ideal. It is well known that if every prime ideal of a commutative ring is a finitely generated ideal then it is a Noetherian ring. Therefore $\Lambda$ is a Noetherian ring. \\
(iv)$\Leftrightarrow$(v): It is well known that a commutative ring is an Artinian ring if and only if it is a Noetherian ring with the Krull dimension zero. \\
(iv)$\Rightarrow$(ii): It is well known that in a Noetherian ring every ideal has a primary decomposition. Thus there exists a finite set $\{\mathfrak{q}_{1},...,\mathfrak{q}_{n}\}$
of primary ideals of $\Lambda$ such that $\bigcap\limits_{k=1}^{n}\mathfrak{q}_{k}=0$. It is easy to see that the radical of every primary ideal of a ring is a prime ideal. By Theorem \ref{Proposition DIII}(iii), every ideal of $\Lambda $ is a radical ideal. Therefore each $\mathfrak{q}_{k}$ is a maximal ideal of $\Lambda$.
\end{proof}

\begin{remark} Let $\Lambda'$ be the set of all $f\in\Lambda$ such that $\Su(f)$ is either finite or cofinite (i.e., its complement is finite). Clearly $1\in\Lambda'$, since $\Su(1)=X$ is cofinite. Let $f,g\in\Lambda'$. Then clearly $fg\in\Lambda'$, because $\Su(fg)=\Su(f)\cap\Su(g)$. If both $\Su(f)$ and $\Su(g)$ are finite then $\Su(f+g)$ is also finite, since $\Su(f+g)\subseteq\Su(f)\cup\Su(g)$. If $\Su(f)$ is finite and $\Su(g)$ is cofinite then $\Su(f+g)$ is cofinite, because $\Su(f+g)^{c}\subseteq\Su(f)\cup\Su(g)^{c}$. Thus in these cases $f+g\in\Lambda'$. But if both $\Su(f)$ and $\Su(g)$ are cofinite then $f+g$ is not necessarily a member of $\Lambda'$. For example, take $X=\mathbb{N}$ and $K_{x}=\mathbb{Z}_{3}$ for all $x$, then the sequence $f:=((-1)^{n})\in\Lambda'$ but $\Su(1+f)$ is neither finite nor cofinite, and so $1+f\notin\Lambda'$. Therefore $\Lambda'$ is not necessarily a subring of $\Lambda$. One can easily see that $\Lambda'$ is a subring of $\Lambda$ if and only if $K_{x}=\mathbb{Z}_{2}$ for all but a finite number of indices $x$. Hence this reduces to study the subring of $\mathcal{P}(X)\simeq\prod\limits_{x\in X}\mathbb{Z}_{2}$ consisting of all subsets of $X$ which are either finite or cofinite. We have already
studied this subring in \cite[\S7]{Tarizadeh} and some interesting results have been obtained.
\end{remark}

Let $T$ be the set of all $f\in\Lambda$ such that $\Su(f)$ is cofinite. Clearly it is a multiplicative set, and we have the following result.

\begin{theorem}\label{Theorem IVD2} The ring $T^{-1}\Lambda$ is canonically isomorphic to $\Lambda/\mathfrak{F}(\Lambda)$.
\end{theorem}

\begin{proof} If $f\in\mathfrak{F}(\Lambda)$ then consider the sequence $g=(g_{x})\in\Lambda$ such that $g_{x}$ is either $0$ or $1$, according as $x\in\Su(f)$ or $x\notin\Su(f)$. Then clearly $fg=0$ and $g\in T$ because $\Su(g)^{c}=\Su(f)$ is finite. Therefore
$\mathfrak{F}(\Lambda)\subseteq\Ker\pi$ where $\pi:\Lambda\rightarrow T^{-1}\Lambda$ is the canonical ring map. (In fact $\mathfrak{F}(\Lambda)=\Ker\pi$). Thus $\pi$ induces a morphism of rings $\phi:\Lambda/\mathfrak{F}(\Lambda)\rightarrow T^{-1}\Lambda$ given by $f+\mathfrak{F}(\Lambda)\mapsto f/1$. The image of each $g\in T$ under the canonical ring map  $\Lambda\rightarrow\Lambda/\mathfrak{F}(\Lambda)$ is invertible, because consider the sequence $g^{\ast}\in\Lambda$ where $g_{x}^{\ast}$ is either $g^{-1}_{x}$ or $0$, according as $x\in\Su(g)$ or $x\notin\Su(g)$ then $1-gg^{\ast}=\sum\limits_{x\in X\setminus\Su(g)}\Delta_{x}\in\mathfrak{F}(\Lambda)$. Thus there exists a (unique) morphism of rings $\psi:T^{-1}\Lambda\rightarrow\Lambda/\mathfrak{F}(\Lambda)$ given by $f/g\mapsto fg^{\ast}+\mathfrak{F}(\Lambda)$. Then clearly $\psi\circ\phi$ is the identity and $(\phi\circ\psi)(f/g)=(f/1)\big((\phi\circ\psi)(g/1)\big)^{-1}=f/g$.
\end{proof}

\begin{corollary} Let $M$ be a maximal ideal of $\Lambda$. Then $M\cap T=\emptyset$ if and only if $\mathfrak{F}(\Lambda)\subseteq M$.
\end{corollary}

\begin{proof} Let $M\cap T=\emptyset$. It suffices to show that $\Delta_{x}\in M$ for all $x\in X$. Suppose there exists some $x\in X$ such that $\Delta_{x}\notin M$. This yields that $1-\Delta_{x}\in M$. But this is a contradiction since $1-\Delta_{x}\in T$. Conversely, if $f\in M\cap T$ then $\sum\limits_{x\in X\setminus\Su(f)}\Delta_{x}\notin M$. Thus there exists some $x\in X\setminus\Su(f)$ such that $\Delta_{x}\notin M$. But this is a contradiction since $\mathfrak{F}(\Lambda)\subseteq M$. \\
As a second proof, by Theorem \ref{Theorem IVD2}, we have $\{\mathfrak{p}\in\Spec(\Lambda):\mathfrak{F}(\Lambda)
\subseteq\mathfrak{p}\}=
\{\mathfrak{p}\in\Spec(\Lambda):\mathfrak{p}\cap T=\emptyset\}$.
\end{proof}

If $K=\mathbb{F}_{q}$ is a finite field then $f=f^{q}$ for all $f\in\Fun(X,K)$, since every element of $K$ satisfies in the equation $x=x^{q}$. This fact also leads us to the following result.

\begin{theorem}\label{Th Del} If $K=\mathbb{F}_{q}$ is a finite field, then the residue fields of the ring $R=\Fun(X,K)$ are canonically isomorphic to $K$.
\end{theorem}

\begin{proof} Let $M$ be a maximal ideal of $R$. The canonical ring map $\pi:K\rightarrow R/M$ is injective. Thus $R/M$ has at least $q$ elements. But every element of $R/M$ satisfies in the equation $x=x^{q}$. If every element of an integral domain $R$ satisfies in the equation $x=x^{q}$, then every non-zero element of $R$ is a root of the polynomial $x^{q-1}-1\in R[x]$, thus $R$ has at most $q$ elements. Hence, $R/M$ has exactly $q$ elements, and so $\pi$ is an isomorphism.
\end{proof}

For a given function $f:X\rightarrow Y$ and an integral domain $K$, then it can be easily seen that $\eta_{Y}\circ f=\Fun(f)^{\ast}\circ\eta_{X}$ where $\eta_{X}:X\rightarrow\Spec\big(\Fun(X,K)\big)$ is the canonical map which is given by $x\mapsto\mathfrak{p}_{x}=\{f\in\Fun(X,K): f(x)=0\}$. We also have the following easy observation.

\begin{proposition}\label{Remark I} For a function $f:X\rightarrow Y$ and a nonzero ring $K$ the following assertions hold. \\
$\mathbf{(i)}$ $f$ is injective if and only if $\Fun(f):\Fun(Y,K)\rightarrow\Fun(X,K)$ is surjective. \\
$\mathbf{(ii)}$ $f$ is surjective if and only if $\Fun(f)$ is injective.
\end{proposition}

\begin{proof}
(i): Assume $f$ is injective and take a map $g:X\rightarrow K$. Consider the map $h:Y\rightarrow K$ given by $h(y)=g(x)$ if $y=f(x)$ for some $x\in X$ otherwise $h(y)=0$. Clearly $h$ is well-defined, since $f$ is injective. Obviously $h\circ f=g$. Hence, $\Fun(f)$ is surjective. Conversely, assume $\Fun(f)$ is surjective. Suppose $f(a)=f(b)$ for some $a,b\in X$. We show that $a=b$. Consider the function $g:X\rightarrow K$ with $g(x)=\delta_{a,x}$ where $\delta_{a,x}$ is the Kronecker delta.
Then there exists a map $h:Y\rightarrow K$ such that $h\circ f=g$. In particular, $1=\delta_{a,a}=h(f(a))=h(f(b))=\delta_{a,b}$. It follows that $a=b$. \\
(ii): If $f$ is surjective then clearly $\Fun(f)$ is injective. Conversely, consider the functions $g,h:Y\rightarrow K$ where $g$ is the constant function  1 and $h(y)=1$ if $y$ is in the image of $f$ otherwise $h(y)=0$. Then clearly $g\circ f=h\circ f$. It follows that $g=h$, because $\Fun(f)$ is injective. This shows that $f$ is surjective.
\end{proof}

\section{A geometric result}

In this section it is shown that every set $X$ can be made into a separated scheme, and this scheme is an affine scheme if and only if $X$ is a finite set. \\

First we define a sheaf of (commutative) rings $\mathscr{O}_{X}$ on the discrete space $X$ by sending each (open) subset $U$ of $X$ into $\prod\limits_{x\in U}K_{x}$. The restriction morphisms of this sheaf are the canonical projections.

\begin{lemma} If $f\in\Lambda$ then the ring $\Lambda_{f}$ is canonically isomorphic to $\mathscr{O}_{X}\big(\Su(f)\big)=\prod\limits_{x\in\Su(f)}K_{x}$.
\end{lemma}

\begin{proof} The image of $f$ under the canonical projection $\pi:\Lambda\rightarrow\mathscr{O}_{X}\big(\Su(f)\big)$
is invertible in $\mathscr{O}_{X}\big(\Su(f)\big)$. Hence, there exists a (unique) morphism of rings $\pi'_{f}:\Lambda_{f}\rightarrow\mathscr{O}_{X}\big(\Su(f)\big)$ such that $\pi'_{f}(g/f^{n})=\pi(g)\pi(f)^{-n}$. The map $\pi_{f}'$ is surjective. If $\pi_{f}'(g/f^{n})=0$ then $g_{x}=0$ for all $x\in\Su(f)$. It follows that $fg=0$.
\end{proof}

Then we define the canonical morphism of ringed spaces: $$(\eta,\eta^{\sharp}):(X,\mathscr{O}_{X})\rightarrow(\Spec(\Lambda),
\mathscr{O}_{\Spec(\Lambda)})$$ as follows. The map $\eta:X\rightarrow\Spec(\Lambda)$ sends each $x\in X$ into $\mathfrak{m}_{x}$. If $f\in\Lambda$ then $\eta^{-1}\big(\D(f)\big)=\Su(f)$. Let $\psi_{\D(f)}$ be the composition of the following canonical isomorphisms of rings: $$\xymatrix{\mathscr{O}_{\Spec (\Lambda)}\big(\D(f)\big)\ar[r]^{\:\:\:\:\:\:\:\:\:\:\:\:\:\:\:\simeq}&
\Lambda_{f}\ar[r]
^{\pi'_{f}\:\:\:\:\:\:\:\:\:\:\:\:}
&\mathscr{O}_{X}\big(\Su(f)\big).}$$ If $\D(g)\subseteq\D(f)$ then by Theorem \ref{Proposition DIII}(iv), $\Su(g)\subseteq\Su(f)$, and the following diagram is commutative:
 $$\xymatrix{
\mathscr{O}_{\Spec(\Lambda)}\big(\D(f)\big)\ar[r]^{\:\:\:\:\psi_{\D(f)}} \ar[d] &\mathscr{O}_{X}\big(\Su(f)\big)
\ar[d]\\\mathscr{O}_{\Spec(\Lambda)}
\big(\D(g)\big)\ar[r]^{\:\:\:\:\psi_{\D(g)}} &\mathscr{O}_{X}\big(\Su(g)\big)}$$
where the vertical arrows are the restriction morphisms. Hence there exists a (unique) morphism of sheaves of rings $\eta^{\sharp}:\mathscr{O}_{\Spec (\Lambda)}\rightarrow\eta_{\ast}\mathscr{O}_{X}$ such that $\eta_{\D(f)}^{\sharp}=\psi_{\D(f)}$ for all $f\in\Lambda$. In fact, $\eta^{\sharp}$ is an isomorphism of sheaves of rings.

\begin{theorem}\label{Theorem set-scheme} If $X$ is a set then the ringed space $(X,\mathscr{O}_{X})$ is a separated scheme. Moreover, it is an affine scheme if and only if $X$ is a finite set.
\end{theorem}

\begin{proof} If $x\in X$ then the ringed space $(U,\mathscr{O}_{U})$ is isomorphic to the affine scheme $\Spec(K_{x})$ where $U=\{x\}$. Hence, $(X,\mathscr{O}_{X})$ is a scheme. If it is an affine scheme then $X$ is a finite set, because every affine scheme is quasi-compact. Conversely, if $X$ is a finite set then the map $\eta:X\rightarrow\Spec(\Lambda)$ between the underlying spaces is a homeomorphism. Therefore $\eta:X\rightarrow\Spec(\Lambda)$ is an isomorphism of schemes. Thus $X$ is an affine scheme. This shows that the affine opens of the scheme $X$ are precisely the finite subsets of $X$. Using this and Proposition \ref{Proposition DII}(ii), then we get that the scheme $X$ is a separated scheme.
\end{proof}

By the above theorem, the stalk $\mathscr{O}_{X,x}$ is canonically isomorphic to the field $K_{x}$. \\

For a fixed field $K$, then the scheme $(X,\mathscr{O}_{X})$ with the structure sheaf defined by $\mathscr{O}_{X}(U)=\Fun(U,K)$ with $U\subseteq X$ has the functorial aspect. In fact, if $f:X\rightarrow Y$ is a function and $U$ is a subset of $Y$ then we define $f_{U}^{\sharp}:=\Fun(f_{U})$ where the function $f_{U}:f^{-1}(U)\rightarrow U$ is induced by $f$. It is easy to see that $(f,f^{\sharp}):(X,\mathscr{O}_{X})\rightarrow(Y,\mathscr{O}_{Y})$ is a morphism of schemes. In fact, the assignments $X\mapsto(X,\mathscr{O}_{X})$ and $f\mapsto(f,f^{\sharp})$
define a (faithful) covariant functor from the category of sets to the category of schemes.


\begin{thebibliography}{10}
\bibitem{Anderson}
D.D. Anderson and J. Kintzinger, Ideals in direct products of commutative rings, Bull. Austral. Math. Soc. 77 (2008) 477-483.
\bibitem{Gilmer-Heinzer}
R. Gilmer and W. Heinzer, On the imbedding of a direct product into a zero-dimensional commutative ring, Proc. Amer. Math. Soc. 106 (1989) 631-637.
\bibitem{Gilmer-Heinzer 2}
R. Gilmer and W. Heinzer, Products of commutative rings and zero-dimensionality, Trans. Amer. Math. Soc. 331 (1992) 663-680.
\bibitem{Gilmer-Heinzer 3}
R. Gilmer and W. Heinzer, Infinite products of zero-dimensional commutative rings, Houston J. Math. 21(2) (1995) 247-259.
\bibitem{Maroscia}
P. Maroscia, Sur  les anneaux de dimension zero, Atti Accad. Naz. Lincei Rend Cl. Sei. Fis. Mat. Natur. 56 (1974) 451-459.
\bibitem{Ronnie et al}
R. Levy, et al., The prime spectrum of an infinite product of copies of $\mathbb{Z}$, Fund. Math. 138 (1991) 155-164.
\bibitem{Abolfazl-Taheri}
A. Tarizadeh and Z. Taheri, Stone type representations and dualities by power set ring, J. Pure Appl. Algebra, 225(11) (2021) 106737.
\bibitem{Tarizadeh}
A. Tarizadeh and M.R. Rezaee, Minimal and maximal spectra as the Stone-\v{C}ech compactification, arXiv.1910.09884, (2019).
\bibitem{Tarizadeh 3}
A. Tarizadeh and N. Shirmohammadi, Tame and wild primes in direct products of commutative rings,
https://doi.org/10.48550/arXiv.2208.08828, (2022).
\end{thebibliography}
\end{document}